\documentstyle[12pt]{article}

\def \E {{\mathcal E}}
\def \H {{\mathcal H}}
\def \Z {{\mathcal Z}}
\def \C {{\bf C}}
\def \P {{\rm Pr}}

\newtheorem{theorem}{Theorem}

\title{Zeros of Gaussian Analytic Functions}

\author{M. Sodin
\thanks{Partially supported by grant number 96-00030 of the United
States - Israel Binational Science Foundation, and by grant number 93/97.1
of the Israel Science Foundation of the Israel Academy of Sciences and 
Humanities.}} 
\date{}
\begin{document}

\maketitle
\section{Introduction}

Throughout this note we shall use the following notation.
Let $G\subseteq \C^1$ be a plane domain and $\{\psi_j(z)\}_{j=1}^N$ 
be a system of $N\le\infty$ analytic functions in $G$. By
$\Psi(z)$
we denote the holomorphic curve in the euclidean space $\C^N$ with
coordinates $\psi_j(z)$. If $N=\infty$, we assume that 
\begin{equation} 
\label{eq:0.1}     
||\Psi(z)|| = 
\sum_j |\psi_j(z)|^2<\infty\,, \qquad z\in G,
\end{equation}
where the series on the RHS converges locally uniformly in $G$. 

Let $\omega_j$ be independent, complex-valued, gaussian random variables
such that
$$
\E\left\{\omega_j\right\}=0, \qquad \rm {and} \qquad
\E\left\{|\omega_j|^2\right\}=1.
$$
We identify the probability space with $\C^N$ equipped with the gaussian
product measure $d\nu$. 

A gaussian analytic function $\psi(z,\omega)$ is defined as 
\begin{equation}
\label {eq:0.2}
\psi(z,\omega)=\sum_j \omega_j \psi_j(z)
\end{equation}
(cf. \cite{Kahane2}).
If $N=\infty$, then according to a theorem by Khintchin and Kolmogorov 
(cf. \cite[Chapter~3, Theorem~2]{Kahane}), 
the series converges locally uniformly in $G$ and almost surely in $\omega$, 
and hence defines an
analytic function in $G$.

Let $n_{\omega}$ be a counting measure of zeros (according to their
multiplicities) of the function $\psi(z,\omega)$. Here we shall be
concerned with three general results on  the random measure
$n_{\omega}$.
The first one is a formula for the average $\E\{n_{\omega}\}$
which is due to Edelman and Kostlan (cf. \cite[Theorem~8.1]{EK}). 
The second, close to Calabi's Rigidity Theorem \cite{Calabi},
loosely speaking says that the average measure $\E\{n_\omega\}$
``almost determines'' the analytic functions $\psi_j(z)$. The
third result, which is due to Offord \cite{Offord1} (cf. \cite{Offord2}
and \cite{Offord3}), is an exponential decrease 
of ``tail probabilities'' of an analytic function having an excess or
deficiency of zeros in a given region.
An important feature of these three results is that they do not need any
assumptions about analytic functions $\psi_j(z)$ and ``dimension'' $N$.
We shall not touch on the more delicate statistics of the
local correlation functions, which was recently of some interest
in mathematical physics (see references at the end of this note). 

By $C$ and $c$ we denote various positive numerical constants which may
vary from line to line.

\section{Edelman-Kostlan formula for the average number of zeros}

\begin{theorem} In the assumptions formulated above,
\begin{equation}
\label{eq:1.1}
\E\{n_\omega\}=\frac{1}{2\pi} \Delta \log||\Psi||\,dm(z)\,,
\end{equation}
where $\Delta$ is a distributional Laplacian and $dm(z)$ is the plane
Lebesgue measure; i.e., the average
$\E\{n_\omega\}$ coincides with the Riesz measure of the subharmonic
function $\log||\Psi(z)||$.
\end{theorem}

For $N<\infty$,
Edelman and Kostlan outlined the proof. In a
special case, the proof appears in \cite[Lemma~3.1]{SZ}. 

\medskip
It is curious to note that since
\begin{equation}
\label{eq:1.1a}
\E\left\{|\psi(z,\omega)|^2\right\}
=\sum_{i,j}\E\{\omega_i\overline{\omega_j}\}\psi_i(z)\overline{\psi_j(z)}
=||\Psi(z)||^2\,,
\end{equation}
equation (\ref{eq:1.1}) can be rewritten
in the form
$$
\E\left\{ \frac{1}{2\pi} \Delta\log|\psi(z,\omega)|^2 \right\}
=\frac{1}{2\pi} \Delta \log \E\left\{|\psi(z,\omega)|^2 \right\}\,.
$$

\medskip\par\noindent{\em Proof of Theorem~1:} 
Let $\phi$ be a test
function with a
compact support in $G$. Then, according to Green's formula,
\begin{equation}
\label{eq:1.2}
\int_G \phi(z) dn_\omega (z) = 
\frac{1}{2\pi} \int (\Delta \phi)(z) \log|\psi(z,\omega)|\,dm(z)\,.
\end{equation}

First, we assume that $N<\infty$. In this case we repeat verbatim the
argument from \cite[Lemma~3.1]{SZ}. We integrate the both sides of
formula (\ref{eq:1.2}) against the gaussian measure $d\nu$ in $\C^N$. 
Changing the order of integration on the RHS and using the notation 
$(a,b)=\sum_j a_jb_j$ for the scalar product in $\C^N$ without complex
conjugation, and $\hat a = a/||a||$, we obtain
\begin{eqnarray}
\int_G \phi(z)  \E\{dn_\omega(z)\} 
&=&
\frac{1}{2\pi}\int_G (\Delta \phi)(z)\,dm(z)\, 
\int_{\C^N} \log |(\Psi(z),\omega)|\, d\nu(\omega) \nonumber
\\ \nonumber \\
&=& 
\frac{1}{2\pi}\int_G (\Delta \phi)(z)\,dm(z)\, 
\Big\{
\log||\Psi(z)||
                     \nonumber \\ \nonumber \\  
&&+ \int_{\C^N} \log ||\omega||\, d\nu(\omega) 
+ \int_{\C^N} \log |(\hat\Psi(z),\hat\omega)|\, d\nu(\omega) 
\Big\}\,. 
            \nonumber
\end{eqnarray}
The second integral in the braces does not depend on $z$. The third
integral in the braces also does not depend on $z$ because of
the rotational
invariance of the gaussian measure $d\nu$. Hence, by the Gauss formula,
after integration against $(\Delta\phi)(z)$, where $\phi$ is compactly
supported in $G$, these two terms vanish, and we obtain
$$
\int_G \phi(z) \E\{dn_\omega(z)\} =
\frac{1}{2\pi}\int_G (\Delta \phi)(z)\log||\Psi(z)|| \,dm(z)\,, 
$$  
which is equivalent to (\ref{eq:1.1}).

Now, consider the case $N=\infty$. Clearly, we no longer have  rotational
invariance of $d\nu$. Instead, we shall use an approximation
argument together with the following fact which will be proved in
Section~4:
for each $z\in G$, and each $\lambda>0$, 
\begin{equation}
\label{eq:1.3}
\P\left(\left\{
\omega:\, \left|\log|\psi(z,\omega)|-\log||\Psi(z)||\,\right|
>\lambda
\right\}\right)\le 3e^{-\lambda}\,.
\end{equation} 
This yields that for every compact subset $K$ of $G$
\begin{equation}
\label{eq:1.4}
\int_{\C^\infty}d\nu(\omega)\,\int_K \log^2|\psi(z,\omega)| \,dm(z)
\le Cm(K)+ 2\int_K \log^2||\Psi(z)||\, dm(z)\,.
\end{equation}

Let $k$ be a positive integer. Denote by $\Psi_k(z)$ a holomorphic curve
in $\C^k$ which consists of the first $k$ components of the curve
$\Psi(z)$, and set $\psi_k(z,\omega)=(\Psi_k(z),\omega)$. For a.a.
$\omega\in\C^\infty$, the sequence of analytic functions 
$z\mapsto \psi_k(z,\omega)$ converges to $\psi(z,\omega)$ locally
uniformly in $G$, and therefore
\begin{equation}
\label{eq:1.5}
\lim_{k\to\infty} \int_K 
\big| \log|\psi(z,\omega)|-\log|\psi_k(z,\omega)| \big|
\,dm(z) = 0\,.  
\end{equation}  
On the other hand, applying estimate (\ref{eq:1.4}) to the functions 
$\psi_k(z,\omega)$, and using the locally uniform convergence of
$\log||\Psi_k(z)||$ to $\log||\Psi(z)||$, we get:
\begin{equation}
\label{eq:1.6}
\int_{\C^\infty}d\nu(\omega)\,\int_K \log^2|\psi_k(z,\omega)| \,dm(z)
\le Cm(K)+ 3\int_K \log^2||\Psi(z)||\, dm(z)\,,
\end{equation}
if $k$ is sufficiently large.

Put
$$
h_k(\omega)=\int_K
\left| \log|\psi_k(z,\omega)|-\log|\psi(z,\omega)|\,
\right|\, dm(z)\,.
$$
Then, by (\ref{eq:1.5}), 
\begin{equation}
\label{eq:1.7}
\lim_{k\to\infty} h_k(\omega) = 0\,, \qquad 
\nu-{\rm a.e}\,. 
\end{equation}
and by (\ref{eq:1.4}) and (\ref{eq:1.6})
\begin{equation}
\label{eq:1.8}
\int_{\C^\infty} h_k^2(\omega)\, d\nu(\omega) \le M(K)\,,
\end{equation}
where the constant $M(K)$ does not depend on $k$. 

Relations (\ref{eq:1.7}) and (\ref{eq:1.8}) yield that
$$
\lim_{k\to\infty} \int_{\C^\infty} h_k(\omega)\, d\nu(\omega) = 0\,,
$$
or
$$
\lim_{k\to\infty} \int_K dm(z) \int_{\C^\infty} 
\left| \log|\psi_k(z,\omega)|-\log|\psi(z,\omega)|\,
\right|\, d\nu(\omega)
= 0\,.
$$
Since $K$ is an arbitrary compact subset of $G$, this yields that the
sequence of subharmonic functions 
$$
z\mapsto \int_{\C^\infty} \log|\psi_k(z,\omega)|\,d\nu(\omega)\,,
\qquad k=1,2,...\,, \quad z\in G\,,
$$
converges in $L^1_{\rm{loc}}(G)$ to the subharmonic function
$$
z\mapsto \int_{\C^\infty} \log|\psi(z,\omega)|\,d\nu(\omega)\,,
\qquad z\in G\,.
$$
By continuity of the Laplacian in $L^1_{\rm{loc}}(G)$, this yields weak
convergence of the corresponding Riesz measures. That is, for
$k\to\infty$, 
the sequence $\E\{n_{\omega,k}\}$
weakly converges to $\E\{n_\omega\}$, 
where $n_{\omega,k}$ is the counting measure of zeros of
$\psi_k(z,\omega)$. 
Besides, the Riesz measures of $\log||\Psi_k(z)||$ also converge weakly to 
the Riesz measures of $\log||\Psi(z)||$. By the first part of the proof,  
$\E\{n_{\omega,k}\}$ coincides with the Riesz measure of
$\log||\Psi_k(z)||$. 
This completes the proof. 
$\Box$ 

\section{Calabi's Rigidity}

In this section, we assume that the functions $\psi_j$ are linearly
independent, that is, if
$$
\sum_j c_j\psi_j(z) \equiv 0\,, \qquad \sum_j |c_j|^2<\infty\,,
$$
then all $c_j$ must vanish. 

We are interested in the following question
\begin{itemize}
\item[$\bullet$] what we are allowed to do with the functions
$\psi_j(z)$ witnout changing the average measure $\E\{n_{\omega}\}$?
\end{itemize}
Surely, we are allowed to multiply all functions $\psi_j$ by the same
analytic function $g$ without zeros in $G$. We can also apply a constant
unitary transformation $U$ to the vector-function $\Psi(z)$ since this
does not change the distribution  of the function $\psi(z,\omega)$. 
Amazingly enough, the
measure $\E\{n_\omega\}$ determines the functions $\psi_j (z)$ (and the
dimension $N$) up to these operations.  

\begin{theorem}
Let $\Psi_1(z)$ and $\Psi_2(z)$ be analytic vector functions in $G$ with
linearly independent components and with the same average measure
$\E\{n_\omega\}$ of zeros of $(\Psi(z),\omega)$. Then, the dimension $N$
of the vectors $\Psi_1$ and $\Psi_2$ is the same, and there exists a
scalar analytic function $g(z)$ without zeros in $G$, and a unitary
transformation $U$ in ${\bf C}^N$ such that
\begin{equation}
\Psi_2(z) = g(z) U \Psi_1 (z)\,.
\label{star}
\end{equation}
\end{theorem}

\medskip\par\noindent{\em Proof:} Having a vector-function $\Psi(z)$ with
linearly independent analytic components $\psi_j(z)$, define the Hilbert
space $\H$ of analytic functions in $G$ with elements
$$
f(z) = \sum_j c_j \psi_j (z)\,, \qquad
||f||_{\H}=\sum_j |c_j|^2\,.
$$
This is a kernel space with the reproducing kernel 
$$
K(z,w)=\sum_j \psi_j(z) \overline{\psi_j (w)}\,,
\qquad f(w)=<f,K(\,.\,,w)>_{\H}\,.
$$
The function $K(z,w)$ is analytic in $z$, anti-analytic in $w$ and is
hermitian-positive.
According to Theorem~1, the measure $\mu=\E\{n_\omega\}$ coincides with
the Riesz measure of the subharmonic function $\log\sqrt{K(z,z)}$.

Now, assume that we have two Hilbert spaces $\H_1$ and $\H_2$ of
analytic functions in $G$ with kernels $K_1$ and $K_2$ such that
$\log\sqrt{K_1(z,z)}$ and $\log\sqrt{K_2(z,z)}$ have the same Riesz
measure, that is, the function $\log K_1(z,z)-\log K_2(z,z)$ is harmonic
in $G$. Then
\begin{equation}
K_2(z,z) = e^{2u(z)} K_1(z,z)\,,
\label{alpha}
\end{equation}
where $u(z)$ is harmonic in $G$. 

Observe that the diagonal 
$\{(\lambda,\lambda):\, \lambda\in G\}$
is a set of uniqueness for functions $K(z,w)$ analytic in $z\in G$ and 
anti-analytic in $w\in G$. Indeed, for a moment set ${\widehat
K}(\lambda)=
K(\lambda,\lambda)$. Then, for $m,n\in {\bf Z_+}$,
$$
\frac{\partial^{m+n}K}{\partial z^m \partial
{\bar w}^n}\Big|_{(\lambda,\lambda)}
=\frac{\partial^{m+n}{\widehat K}(\lambda)}{\partial \lambda^m \partial
{\bar \lambda}^n}\,.
$$
This follows from the formal rules of taking $\partial$ and
$\bar\partial$ derivatives.
(It is enough to check this for monomials of the form $z^k{\bar w}^l$,
$k,l\in{\bf Z_+}$.) Therefore, knowing  ${\widehat K}(\lambda)$, we may
use
the
Taylor formula and hence extend the kernel $K$ from the diagonal to its
neighbourhood.

This observation shows that
\begin{equation}
K_2(z,w)=g(z)\overline{g(w)} K_1(z,w)\,,
\label{beta}
\end{equation}
where $g(z)$ is an analytic function in $G$ without zeros, and
$\log|g|=u$. Indeed, the RHSs of (\ref{alpha}) and
(\ref{beta}) coincide on the diagonal. Since the both kernels $K_1$ and
$K_2$ are single-valued functions in $z$ and $w$, the function $g$ is
also single-valued.

Having at hands the kernel $K$, one recovers the space $\H$ as the closure
of finite linear combinations 
$$
\zeta(z) =\sum_i \zeta_iK(z,w_i)\,, \qquad \{\zeta_i\}\subset\C\,,
\quad \{w_i\}\subset G\,,
$$
in the norm
$$
||\zeta||_{\H}^2 = \sum_{i,j} \zeta_i\overline{\zeta_j} K(w_i,w_j)\,.
$$
Then equation (\ref{beta}) yields that
\begin{equation}
\H_2 = \left\{f=gh:\, h\in \H_1 \right\}\,, \qquad
<f_1,f_2>_{\H_2}=<h_1,h_2>_{\H_1}\,.
\label{gamma}
\end{equation}
By construction, the components $\{\psi_{1,j}\}$ of $\Psi_1$ form an
orthonormal basis in $\H_1$, while the components of $\Psi_2$ form an
orthonormal basis in $\H_2$. Therefore, relations (\ref{gamma}) are
equivalent to (\ref{star}). This completes the proof. $\Box$

\medskip An equivalent formulation can be made in the framework of
differential geometry.    
As above, consider a holomorphic curve $\Psi (z)$ in $\C^N$ with
components $\psi_j(z)$. The curve gives a holomorphic embedding of $G$
equipped with a Riemannian metric into euclidean
$N$-dimensional space. The argument given above shows that knowing the
norm of the curve
$$
||\Psi (z)||^2 = \sum_j |\psi_j (z)|^2
$$
in a subdomain $G'\subseteq G$, we may recover the dimension $N$ and the
whole curve $\Psi (z)$ up to a constant unitary transformation of the
space
$\C^N$. This is a special case of Calabi's rigidity theorem
\cite{Calabi} (cf. Problem~207 and reference to its solution in Part~IV of
the revised  edition of the P\'olya and
Szeg\"o problem book \cite[pp.34,211]{PS}). 

\medskip The Riesz measure $\mu_K$ of the subharmonic function
$\log\sqrt{K(z,z)}$ is an important invariant of the kernel Hilbert space
$\H$ of analytic functions \cite{CD}. In this setting, a fact equivalent
to Theorem~2 was also found by Nikolskii \cite[Theorem~2.4 and
Lemma~2.5]{nik}). 

\medskip Let $K(z,w)$ be an arbitrary hermitian-positive
kernel, analytic in $z\in G$, and anti-analytic in $w\in G$. As above, let
$\mu_K$ be the Riesz measure of $\log\sqrt{K(z,z)}$. According to
Theorem~2, the class of such measures possesses a strong uniqueness
property: if $\mu_{K_1}=\mu_{K_2}$ in a small disk $D$ in $G$, then
$\mu_{K_1}-\mu_{K_2}$ is a divisor of a meromorphic function in $G$, and
the functions $K_1$ and $K_2$ are related to each other as in
(\ref{beta}) with a meromorphic function $g$ without zeros and poles in
$D$. It would be interesting to find a {\em quantative} version
of this statement when the difference $\mu_{K_1}-\mu_{K_2}$ is small in
$D$. 

In \cite{Calabi1}, Calabi completely decribed the class
of all measures $\mu_K$, see also Lawson's paper \cite{Lawson}. It
seems, this description deserves a better understanding.

\section{Offord's estimate for large deviations}

\begin{theorem} Let $\mu$ be the Riesz measure of the subharmonic
function $\log ||\Psi(z)||$, and let $\phi\in C_0^\infty (G)$ be an
arbitrary
test function with compact support in $G$. Then, for every $\lambda>0$, 
\begin{equation}
\label{eq:3.1}
\P\left(\left\{\omega:\, 
\left| 
\int_G \phi (dn_{\omega} - \E\{dn_{\omega}\})
\right|\ge \lambda \right\}\right)
\le 3e^{-2\pi\lambda/||\Delta\phi||_{L^1}}\,.
\end{equation}
\end{theorem}

We start with a lemma which is a key ingredient in Offord's approach.

\medskip\par\noindent{\bf Lemma} {\em Let $\Z (\omega)$ be a
complex gaussian random variable with zero average and variance 
$\sigma^2$.
Then, for every measurable subset $E$,}
\begin{equation}
\label{eq:3.2}
\left|\int_E \log|\Z(\omega)|\,d\nu(\omega) - \nu(E)\log\sigma 
\right|
\le \nu(E)\left[\log\frac{1}{\nu(E)}+\frac{1}{4}\right]\,.
\end{equation}

\smallskip\par\noindent{\em Proof of Lemma:} WLOG, we may assume
that $\sigma=1$, otherwise we replace $\Z$ by $\Z/\sigma$. The upper
bound follows at once from Jensen's inequality:
\begin{eqnarray}
\int_E \log|\Z|\,d\nu &=& 
\frac{\nu(E)}{2}\,\frac{1}{\nu(E)}\, \int_E \log|\Z|^2\,d\nu
\nonumber \\ 
&\le& 
\frac{\nu(E)}{2}\,\log\left(\frac{1}{\nu(E)}\, \int_E |\Z|^2\,d\nu
\right) \le \frac{\nu(E)}{2}\,\log\frac{1}{\nu(E)}\,.
\nonumber
\end{eqnarray}

Now, let us prove the lower bound in (\ref{eq:3.2}). We have
\begin{eqnarray}
\int_E \log|\Z|\,d\nu &\ge& 
- \int_E \log^-|\Z|\,d\nu 
\nonumber \\ \nonumber \\
&=&- \int_{E,\, |\Z|\le \nu(E)} \log^-|\Z|\,d\nu 
- \int_{E,\, |\Z|> \nu(E)} \log^-|\Z|\,d\nu 
\nonumber \\ \nonumber \\
&\ge& -\int_0^{\nu(E)} \frac{\P(|\Z|\le s)}{s}\,ds
-\nu(E)\log\frac{1}{\nu(E)}\,. \nonumber
\end{eqnarray}

Since $\Z$ is gaussian with variance one, we can easily estimate the
first integral in the RHS:
\begin{eqnarray}
\label{eq:concentration}
\int_0^{\nu(E)} \frac{\P(|\Z|\le s)}{s}\,ds
&=&\int_0^{\nu(E)} \frac{1-e^{-s^2/2}}{s}\,ds \nonumber \\ \nonumber \\ 
&\le& \frac{1}{2}\int_0^{\nu(E)} s\,ds = \frac{\nu^2(E)}{4} \le
\frac{\nu(E)}{4}\,.
\end{eqnarray}
This completes the proof. $\Box$. 

\medskip Now we turn to the

\smallskip\par\noindent{\em Proof of Theorem~3:} Fix $\lambda>0$, and 
define two sets:
$$
A_+=\left\{\omega:\, \int_G \phi\, (dn_\omega-d\mu) \ge \lambda
\right\}
$$
and
$$
A_-=\left\{\omega:\, \int_G \phi\, (dn_\omega-d\mu) \le -\lambda
\right\}\,.
$$

Using the lemma, we can easily estimate $\nu(A_{\pm})$ from the above.
First,
consider $A_+$. Using Green's formula, we obtain
$$
\int_G \phi(z)\, (dn_\omega(z)-d\mu(z))
=
\frac{1}{2\pi}\,
\int_G (\Delta\phi)(z)\left\{ \log|\psi(z,\omega)|-\log
||\Psi(z)|| 
\right\}\, dm(z)\,.
$$
Therefore,
\begin{eqnarray}
\lambda\nu(A_+) &\le& 
\int_{A_+} d\nu(\omega) \, \int_G \phi(z)\, (dn_\omega(z)-d\mu(z))
\nonumber \\ 
&=&
\frac{1}{2\pi}\,
\int_G (\Delta\phi)(z)\, 
\left\{
\int_{A_+}\log|\psi(z,\omega)|\,d\nu(\omega)-\nu(A_+)\log
||\Psi(z)|| 
\right\}\, dm(z)\,. \nonumber
\end{eqnarray}
Then applying the lemma to the random variable
$\Z(\omega)=\psi(\,.\,,\omega)$ and to the set $A_+$ (and recalling
that $||\Psi(z)||$ is a variance of the gaussian random variable
$\psi(z,\omega)$, cf. (\ref{eq:1.1a})), we proceed further:
$$
\le \frac{1}{2\pi}\,
\int_G |(\Delta\phi((z)| \nu(A_+)\left[\log\frac{1}{\nu(A_+)} +
\frac{1}{4}\right] \,dm(z)
$$
$$
=\frac{1}{2\pi}\,
\nu(A_+)\left[\log\frac{1}{\nu(A_+)} +
\frac{1}{4}\right]\,||\Delta\phi||_{L^1}\,. 
$$
That is,
$$
2\pi\lambda\nu(A_+) \le  
\nu(A_+)\left[\log\frac{1}{\nu(A_+)}+\frac{1}{4}\right]\,||\Delta\phi||_{L^1}\,,
$$
or
$$
\nu(A_+)\le e^{-2\pi\lambda/||\Delta\phi||_{L^1} + 1/4}\,.
$$
The same estimate holds for the set $A_-$, and this completes the proof.
$\Box$

\medskip
Observe that the lemma also yields  estimate (\ref{eq:1.3}) used in
Section~2. For this, we apply the lemma to the same random variable
$\Z(\omega)=\log|\psi(\,.\,,\omega)|$ and to the sets
$$
A_+=\left\{\omega:\, \log|\Z(\omega)|-\log||\Psi(z)|| >\lambda 
\right\}\,, 
$$
and 
$$ 
A_-=\left\{\omega:\, \log|\Z(\omega)|-\log||\Psi(z)|| <-\lambda 
\right\}\,.
$$  

\medskip We also observe, that applying the theorem with $\lambda=\int
\phi d\mu$, we obtain
$$
\P\left( \left\{\omega:\,\int \phi\, dn_{\omega}=0 \right\}\right)
\le 3\exp\left[-2\pi\,\int\phi d\mu\Big/||\Delta\phi||_1\right]\,.
$$ 

In turn, this yields an upper bound for the ``hole probability'' 
$\P(\{\omega:\, n_{\omega}(D_R)=0 \} )$.
Let $D_R\subset G$ be a disk of radius $R$, and let $D_r$,
$r<R$, be a concentric disk of a smaller radius $r$. 
Then we may choose
a non-negative test-function 
$\phi (z) = \Phi(|z|)$ equal $1$ on $D_r$, vanishing outside $D_R$, and
such that
$$
||\Delta \phi||_{L^1}
= \int_r^R \left(t|\Phi''(t)| + |\Phi'(t)| \right)\, dt <
C\frac{R+r}{R-r}\,.  
$$
We obtain 

\medskip\par\noindent{\bf Corollary 1.} 
$$
\P\left(\left\{\omega:\, n_{\omega}(D_R)=0 \right\} \right)
\le 3\exp\left[-c\mu(D_r)\frac{R-r}{R+r} \right]\,,
\qquad 0<r<R\,.
$$

\medskip
For example, if $\psi_j(z)$ are analytic in the unit disk $D$, and
$$
p=\P\left(\left\{\omega:\, n_{\omega}(D)=0 \right\} \right) >0\,,
$$
then we get a dimensionless bound
$$
\mu (D_r) \le C\left(\log\frac{3}{p}\right) \cdot \frac{1}{1-r}\,, 
\qquad 0<r<1\,.
$$

\medskip There is a certain resemblance between Offord's theorem and known
results of the value-distribution theory due to Littlewood
\cite{Littlewood} and Ahlfors \cite{Ahlfors}. The reader may also have a
look at Favorov's papers \cite{Favorov1}, \cite{Favorov2}, where 
infinite-dimensional versions of these results were obtained.  
We shall not pursue this matter here.

\section{Remarks}

\medskip\par\noindent{\em 1. Functions of several complex variables.}
The facts about zeros of gaussian analytic functions discussed above,
after simple modifications, remain valid for gaussian analytic functions
of several complex variables. 

\medskip\par\noindent{\em 2. Holomorphic sections of line
bundles.}
There is no need to assume that the analytic
functions $\psi_j(z)$ are single-valued. Instead,
one can deal, say, with zeros of character-automorphic functions or
differential forms. More generally, one can work with zeros of gaussian
holomorphic sections of a holomorphic line bundle.  Since the proofs
given above are ``local'', they carry over the case of  holomorphic
sections as well.  

\medskip\par\noindent{\em 3. Non-gaussian distributions.}
In the proof of Theorem~1, one can replace gaussian distributions by any
rotational invariant distribution in $\C$ which does not charge the
origin. 

The Offord estimate is much more robust. The
gaussian distribution was used only once during the estimate of
concentration (\ref{eq:concentration}). 
For possible generalization of the lemma from Section~4 to more general 
classes of random variables see Offord's papers.
A somewhat similar problem was treated by Favorov, Gorin 
and Ullrich, see \cite{Favorov3}, and references given there.

\medskip\par\noindent Added in proof:

\smallskip\par\noindent
Recently, A. Volberg and the author found another generalization of this
lemma. Let $P$ be a polynomial in $\bf R^n$  of degree
$d$, and let $d\nu$
be a logarithmically concave probability measure in ${\bf R}^n$ (that is, 
$d\nu=e^{-g(x)}dx$, where the set $K=\{g<+\infty\}$ and the function
$g|_K$ are convex). Then, for an arbitrary measurable set $E\subset 
{\bf R}^n$,
$$
\left| \frac{1}{\nu(E)}\int_E \log|P|d\nu - 
\int_{{\bf R}^n} \log|P| d\nu \right|
\le 2d\log\frac{C}{\nu(E)}\,.
$$  
This allows us to extend Theorem~3 to families of analytic functions which
depend polynomially on real parameters.

\bigskip
\centerline{\bf Acknowledgements} 

Zeev Rudnick showed me a
preprint version of the paper by Shiffman and Zelditch \cite{SZ} and thus
stimulated my interest in the subject. I thank Efim Gluskin, Leonid
Pastur, Leonid Polterovich, Zeev Rudnick and Peter Yuditskii for useful
discussions and comments.

\medskip\par\noindent School of Math. Sciences, 
\newline Tel Aviv University, \newline
Ramat Aviv 69978, 
\newline Israel 
\smallskip\par\noindent sodin@math.tau.ac.il
\end{document}